\documentclass[11pt]{amsart}
\usepackage{standard}
\newcommand{\Slag}{S_\lag}

\newcommand{\Psih}{\Psi_h}
\newcommand{\lag}{\mathcal{L}}
\newcommand{\md}{\mathcal{M}}
\newcommand{\Ibar}{\overline{I}}

\newcommand{\ombar}{\overline{\omega}}
\newcommand{\orbit}{\rho}

\title{Spreading of Lagrangian regularity on rational invariant tori}
\author{Jared Wunsch}
\date{\today}

\begin{document}
\begin{abstract}
Let $P_h$ be a self-adjoint semiclassical pseudodifferential operator on a
manifold $M$ such that the bicharacteristic flow of the principal symbol on
$T^*M$ is completely integrable and the subprincipal symbol of $P_h$
vanishes.  Consider a semiclassical family of eigenfunctions, or, more
generally, quasimodes $u_h$ of $P_h.$ We show that on a nondegenerate
rational invariant torus, Lagrangian regularity of $u_h$ (regularity under
test operators characteristic on the torus) propagates both along
bicharacteristics, and also in an additional ``diffractive'' manner.  In
particular, in addition to propagating along null bicharacteristics,
regularity fills in the interiors of small annular tubes of
bicharacteristics.
\end{abstract}

\maketitle

\section{Introduction}
It is a well-known fact of semiclassical microlocal analysis, that the
analogue of H\"ormander's theorem on propagation of singularities for
operators of real principal type \cite{Hormander:Propagation} holds for the
semiclassical wavefront set (also known as ``frequency set''): it
propagates along null bicharacteristics of operators with real principal
symbol \cite{Evans-Zworski,Martinez}.  Given a Lagrangian submanifold
$\lag$ of $T^*M,$ we may introduce a finer notion of regularity, the local
\emph{Lagrangian} regularity along $\lag.$ We show here that on rational
invariant tori in integrable systems, local Lagrangian regularity not only
propagates along bicharacteristics, but spreads in additional ways as well.

Let $P_h$ be a semiclassical pseudodifferential operator on a manifold $M,$
with real principal symbol $p$ (this is automatic if $P$ is self-adjoint).
Assume that the bicharacteristic flow of $p$ is completely integrable.  (In
fact we only need to assume integrability \emph{locally,} near one
invariant torus.)  Let $u_h$ be a family of quasimodes of $P_h,$ i.e.\
assume that $\norm{(P_h-\lambda)u_h}_{L^2} =O(h^N)$ for some $N\in \NN,$ as
$h\downarrow 0$ either through a discrete sequence or continuously. (Note
that this certainly includes the possibility of letting $u_h$ be a sequence
of actual eigenfunctions).  Let $\lag$ be an invariant torus in the
characteristic set $\{p=\lambda\}.$ Then the bicharacteristic flow is by
definition tangent to $\lag,$ and we show (even in the absence of the
integrability hypothesis) that Lagrangian regularity propagates along
bicharacteristics---this is Theorem~A below.  If a single trajectory is
dense in $\lag,$ then this is the whole story for propagation, as the set
on which Lagrangian regularity holds is open, hence the whole torus either
enjoys Lagrangian regularity or none of it does.  At the opposite extreme,
if $\lag$ is a torus on which all frequencis of the motion are rationally
related, we may ask the finer question: what subsets of the space of all
orbits may carry Lagrangian regularity?  The answer (assuming a
nondegeneracy condition holds) turns out to be somewhat constrained: given
a single orbit, Lagrangian regularity along a small tube around it implies
Lagrangian regularity along the orbit itself.  This is our Theorem~B. (In
the special case of two-dimensional tori, we can go further: again, either
the whole torus enjoys Lagrangian regularity or no points on it do.)  The
order of regularity up to which our result holds is constrained by the
order of the quasimode.  We speculate that a finer theorem may be
obtainable by more authentically ``second-microlocal'' methods.

\begin{example}
As a simple example of our main result, Theorem~B, we consider the case
$M=S_x^1 \times S_y^1,$ $P_h = h^2 \Lap =-h^2(\pa^2/\pa x^2+ \pa^2/ \pa y^2);$
we consider Lagrangian regularity on the Lagrangian torus $\lag =
\{\xi=0,\eta =1\}$ for quasimodes satisfying
$$
(h^2\Lap -1) u_h \in h^{k+1} L^2(S^1 \times S^1).
$$ Lagrangian regularity on this particular $\lag$ is special in that we
may test for it using powers of the \emph{differential} operator
$D_x=i^{-1} (\pa/\pa x).$ The theorem tells us the following in this case:
let $\Upsilon(x)$ be a smooth cutoff function supported on $\{\abs{x} \in
[\ep,3\ep]\}$ and nonzero at $\pm 2 \ep.$ Let $\phi$ be another cutoff,
nonzero at the origin and supported in $[-2\ep,2\ep].$ If, for all $k'\leq
k,$ we have
$$
\norm{D_x^{k'}(\Upsilon(x) u_h)}\leq C<\infty,
$$
then  for all $k'\leq k,$
$$
\norm{D_x^{k'}(\phi(x) u_h)}\leq \tilde{C}<\infty,
$$ i.e.\ the $D_x^k$ regularity fills in the ``hole'' in the support
of $\Upsilon.$ In this special case, the result can be proved directly by
employing a positive commutator argument using only differential operators;
the positive commutator will arise from the usual commutant $h^{-1} xD_x.$
\end{example}

A less trivial example, that of the spherical pendulum, is discussed in
\S\ref{section:flow} below.

The methods of proof (and the idea of the paper) arose from work of
Burq-Zworski \cite{Burq-Zworski1,Burq-Zworski2} and a subsequent refinement
by Burq-Hassell-Wunsch \cite{BHW1} on the spreading of $L^2$ mass for
quasimodes on the Bunimovich stadium.  The central argument here is a
generalization of the methods used to prove that a quasimode cannot
concentrate too heavily in the interior of the rectangular part of the
stadium (which is essentially the example discussed above on $M=S^1 \times
S^1$).

We remark that our hypotheses in this paper are quite far from those in the
study of ``quantum integrable systems'' where one examines eigenfunctions
of a system of $n$ commuting operators on an $n$-manifold.  For instance,
if we take $P_h = h^2 \Lap +h^2 V$ on the torus, with $V$ a real valued,
smooth bump-function, then the operator $P_h$ satisfies the hypotheses of
our Theorems A and B, and yet there does not exist a system of $n-1$ other
operators commuting with $P_h,$ with independent symbols.  Moreover, even
in the completely integrable case, given that we study eigenfunctions of a
single operator, it may be possible to use the degeneracy of the system to
construct non- or partially-Lagrangian quasimodes.  Little seems to be
known in this direction.

The author is grateful to Andr\'as Vasy for helpful discussions on
Lagrangian regularity, and to Clark Robinson for introducing him to
isoenergetic nondegeneracy.  He has also benefitted greatly from comments
on an earlier version of the manuscript by Maciej Zworski, and by an
anonymous referee.  This work was supported in part by NSF grant
DMS-0401323.

\section{Lagrangian Regularity}
We begin by setting some notation and recalling some concepts of
semiclassical analysis.  For detailed background on this subject, we refer
the reader to \cite{Evans-Zworski,Martinez}.

Let $M^n$ be a smooth manifold and fix $\lag \subset T^*M$ a Lagrangian
submanifold.  Throughout the rest of the paper, we assume\footnote{We may
just as well assume that $h \downarrow 0$ through a discrete sequence; this
will make no difference in what follows.} $u_h \in L^2(M;\Omega^{1/2}),$
with $h \in (0, h_0);$ here $\Omega^{1/2}$ denotes the bundle of
\emph{half-densities} on $M,$ i.e.\ the square root of the
density bundle $\abs{\Lambda^n M}.$
 We will in future, however, suppress the
half-density nature of $u_h$ as well as its $h$-dependence, writing simply
$u \in L^2(M);$ similarly, all operators will tacitly be semiclassical
families of operators, operating on half-densities.  The hypothesis that
our operators act on half-densities ensures that if $A =\Op_h(a)$ with
$a(x,\xi;h) \sim a_0(x,\xi) + h a_1(x,\xi) +\dots,$ the terms $a_0$
(principal symbol) and $a_1$ (subprincipal symbol) are both invariantly
defined as functions on $T^*M$ (see \cite{Evans-Zworski}).

Furthermore, we will deal with an operator $P$ rather than $P-\lambda,$
absorbing the constant term into the definition of the operator.

We begin by defining a notion of Lagrangian regularity of a family of
functions along $\lag,$ following the treatment of the ``homogeneous'' case
in \cite{Hormander:book4}.
\begin{definition}
Let $\md$ denote the module (over
$\Psih(M)$) of semiclassical pseudodifferential operators with symbols
vanishing on $\lag.$

Let $q \in \lag,$ $k \in \NN,$ and $u \in L^2(M).$ We say that $u$ has
Lagrangian regularity of order $k$ at $q,$ and write $q \in \Slag^k(u),$ if
and only if there is a neighborhood $U$ of $q$ in $T^*M$ such that for all
$k'=0,1,\dots,k$ and all $A_1,\dots, A_{k'} \in \md$ with $\WF' A_j \subset
U,$ $h^{-k'} A_1\cdots A_{k'} u \in L^2(M).$
\end{definition}

\begin{proposition}\label{prop:basis}
Fix $q \in \lag,$ and let $A_i$ ($i=1,\dots,n$) be a
collection of elements of $\md$ with $d\sigma(A_i)$ spanning $N^*_q\lag.$
We have
$$ q \in \Slag(u) \Longleftrightarrow h^{-k'} A_{i_1}\cdots A_{i_{k'}} u
\in L^2\quad \forall (i_1,\dots,i_{k'}) \in \{1,\dots,n\}^{k'},\ k'\leq k.$$
\end{proposition}
\begin{proof}
We begin with the case $k=1.$
Given any $B$ characteristic on $\lag$ and microsupported sufficiently
close to $q,$ we may factor $\sigma(B)=\sum c_i \sigma(A_i)$ by Taylor's
theorem.  Thus, letting $C_i$ be operators with symbol $c_i,$ we obtain
$$
h^{-1} B u= \sum h^{-1} C_i A_i u + R u
$$
for some semiclassical operator $R,$ hence we obtain the desired estimate
on $h^{-1} Bu$ since $R$ is uniformly (in $h$) $L^2$-bounded.

More generally, if $B_{\alpha_1},\dots, B_{i_k}$ is a $k$-tuple of operators
characteristic on $\lag,$ we have
$$
h^{-k} B_{i_1}\cdots B_{i_k} u= h^{-k}\prod_{j=1}^k (C_{i_j}
A_{i_j}+hR_{i_j}) u;
$$
We then obtain the desired estimate inductively, using the fact that each
commutator of the form $[C,A]$ or $[R,A]$ produces a further factor of $h.$
\end{proof}

We note that it follows from the work of Alexandrova \cite{Alexandrova}
that $S_\lag=\lag$ if and only if we can actually write $u$ in the form of
an oscillatory integral
$$
\int a(x,\theta,h) e^{i \phi(x,\theta)/h} \, d\theta
$$ plus a term with semiclassical wavefront set away from $\lag;$ here
$\phi$ is a phase function parametrizing the Lagrangian $\lag$ in the sense
introduced by H\"ormander.  This is the semiclassical analog of a central
result in the H\"ormander-Melrose theory of conic Lagrangian distributions
\cite[Chapter 25]{Hormander:book4}.

We now observe that the analogue of H\"ormander's theorem on propagation of
singularities for operators of real principal type is easy to prove in our
setting. 
\begin{theorema}\label{thm:propagation}
Let $P\in \Psi_h(M)$ have real principal symbol $p.$ Let $\lag\subset
\{p=0\}$ be a Lagrangian submanifold of $T^*M.$ Then $P u\in h^{k+1}
L^2(M)$ implies that $\Slag^k(u)$ is invariant under the Hamilton flow of
$p.$
\end{theorema}

The author is grateful to M.~Zworski for suggesting the following brief proof.
\begin{proof}[Sketch]
By \cite[Theorem~21.1.6]{Hormander:book3}, there is a local
symplectomorphism taking $p$ to $\xi_1$ and $\lag$ to $\lag_0 \equiv
\{\xi=0\}.$ Following the development in \cite{Alexandrova}, we may
quantize this to a semiclassical FIO that conjugates $P$ to $h D_{x^1}$
modulo $O(h^\infty)$ (cf.\ \cite[Theorem~26.1.3]{Hormander:book4} in the
non-semiclassical setting).  Lagrangian regularity along $\lag_0$ is
iterated regularity under $h^{-1} (hD_{x^i}),$ i.e.\ is just classical
Sobolev regularity, uniform in $h.$ The theorem thus reduces to the
statement that Sobolev regularity for solutions to $D_{x^1} u \in
h^{k}L^2(M)$ propagates along the lines $(x^1\in \RR,x' =\text{const}),$ which
is easily verified.
\end{proof}

\section{Integrable flow}\label{section:flow}

We continue to assume that $P \in \Psi_h(M)$ has real principal symbol.  We
now further assume that $p=\sigma(P)$ has \emph{completely integrable}
bicharacteristic flow, i.e.\ that there exist functions $f_2,\dots,f_n$ on
$T^*M,$ Poisson commuting with $p$ and with each other, and with
$dp,df_2,\dots,df_n$ pointwise linearly independent. We again emphasize
that we in fact only require the $f_i$'s to exist in some open subset of
interest in $T^*M.$ Let $\Sigma$ denote the characteristic set in $T^*M.$
Let $(I_1,\dots, I_n,\theta_1,\dots, \theta_n)$ be action-angle variables
and let $\omega_i = \partial p/\partial I_i$ be the frequencies.  We also
let $\omega_{ij} = \pa^2 I/\pa I_i \pa I_j.$ (We refer the reader to
\cite{Arnold} for an account of the theory of integrable systems, and in
particular for a treatment of action-angle variables.)

Let $\lag \subset \Sigma$ be a \emph{rational} invariant torus,
i.e.\ one on which $\omega_i/\omega_j \in \QQ$ for all $i,j=1,\dots,n.$ We
further assume that $\lag$ is \emph{nondegenerate} in the following sense:
we assume that the matrix
\begin{equation}\label{bigmatrix}
\begin{pmatrix}
\omega_{11}& \dots & \omega_{1n} & \omega_1 \\
\vdots & \ddots & \vdots & \vdots\\
\omega_{n1} & \dots & \omega_{nn} & \omega_n\\
\omega_{1} & \dots & \omega_{n} & 0\\
\end{pmatrix}
\end{equation}
is invertible on $\lag.$ This is precisely the condition of
\emph{isoenergetic nondegeneracy} often used in KAM theory (see
\cite{Arnold}, Appendix 8D).  It is easy to verify that the condition is
equivalent to the condition that the map from the energy surface to the
projectivization of the frequencies
$$
\{p=0\} \ni I \mapsto [\omega_1(I):\dots: \omega_n(I)] \in \mathbb{RP}^n
$$
be a local diffeomorphism.

For later convenience, we introduce special notation for the frequencies
 and their derivatives on $\lag:$ we let
$$\ombar_i =\omega|_{\lag},\quad \ombar_{ij}
=\omega_{ij}|_{\lag}.$$

On $\lag,$ we of course only know from Theorem~A that $\Slag^k(u)$ is a union
of orbits of $H_p,$ which, being rational, are not dense in $\lag.$
There are, however, further constraints on $\Slag^k(u).$

\begin{definition}
An \emph{annular neighborhood} of a closed orbit $\orbit$ is an open set
$U=V\backslash K\subset \lag$ such that $\orbit \subset K \subset V$ with
$K$ compact and $V$ open in $\lag.$
\end{definition}

We can now state our main result.
\begin{theoremb}
Suppose $P u=f \in h^{k+1} L^2.$ Let $\orbit$ be a null bicharacteristic
for $p$ on the rational invariant torus $\lag.$ If a small enough annular
neighborhood of $\orbit$ is in $\Slag^k (u),$ then so is $\orbit.$

The meaning of  ``small enough'' depends only on the $\ombar_i$'s.

If $n=2$ then either $\Slag (u) = \lag$ or $\Slag(u) = \emptyset.$
\end{theoremb}
Thus, conormal regularity propagates ``diffusively'' to fill in annular
neighborhoods.

\begin{example}
Horozov \cite{Horozov1,Horozov2} has studied the \emph{spherical pendulum},
i.e.\ the system on $T^*S^2$ with Hamiltonian $h=(1/2)\abs{\xi}^2+x_3$ on
$T^*S^2$ (with $x_3$ one of the Euclidean coordinates on $S^2
\subset\RR^3$).  Integrals of motion are $h$ and $p_\theta,$ the angular
momentum.  Horozov showed that when $h\in (-1,1] \cup
[7/\sqrt{17},\infty),$ all values of $p_\theta$ lead to isoenergetically
nondegenerate invariant tori, while for $h \in (1,7/\sqrt{17}),$ there are
exactly two values of $p_\theta$ for which isoenergetic nondegeneracy
fails.  Thus our results show that if we consider quasimodes for the
operator
$$
P_h = (1/2) h^2\Lap_{S^2} +x_3
$$ then for any torus $\lag$ not associated to the one of the exceptional
pairs of $(h,p_\theta)$ identified by Horozov, either $S_\lag=\lag$ or
$S_\lag=\emptyset.$
\end{example}

\begin{example}
We now illustrate with an example the necessity of the isoenergetic
nondegeneracy condition.  As in the introduction, let $M=S^1 \times S^1,$
but now let $P=hD_x;$ it is easy to verify that \emph{no} Lagrangian torus
is isoenergetically nondegenerate in this case.  Let $\lag =
\{\xi=\eta=0\},$ the zero-section of $T^*M.$ Lagrangian regularity in this
setting is, as noted above, just Sobolev regularity, uniform in $h.$

Let $\psi(y)$ be a bump function supported near $y=0.$ Then $$u(x,y) =
e^{i\psi(y)/\sqrt{h}}$$ has wavefront set only in $\lag.$ It is manifestly
Lagrangian on the complement of $\supp\psi,$ which forms an annular
neighborhood of the orbit $\{x\in S^1,y =0, \xi=\eta=0\}\subset \supp
\psi.$ It is not Lagrangian, however, on $\supp \psi,$ as it lacks iterated
regularity under $h^{-1}(h D_y).$

\end{example}
\section{Symbol Construction}
By shifting coordinates, we may assume that $\orbit$ is the orbit passing
through $\{\theta=0\}.$

For each $i,j$ let
$$ \gamma_{ij}= \widetilde{\min}_{k,l \in \ZZ} ((\theta_i+2\pi k)
\ombar_j-(\theta_j+2\pi l) \ombar_i),
$$ where $\widetilde{\min}$ denotes the value with the smallest norm, i.e.\
may be positive or negative.  Each $\gamma_{ij}$ then takes values in an
interval determined by $\ombar_i,$ $\ombar_j,$ and is smooth where it takes
on values in the interior of the interval.  (If $\ombar_i=p/q$ and
$\ombar_j=p'/q'$ then $\gamma_{ij}$ takes values in $[-\pi a, \pi a]$ where
$a=\gcd(qp',pq')/q q'.$) The ``small enough'' condition in the statement of
Theorem~B is just the following: each $\gamma_{ij}$ should be smooth on the
annular neighborhood of $\orbit$ where we assume Lagrangian regularity.

Note that $\gamma_{ij}(\theta)=0$ for all $i,j$ exactly when there exists
$\tilde{\theta} \in \RR^n,$ equivalent to $\theta$ modulo $2\pi\ZZ^n,$ such
that $[\tilde\theta_1: \dots \tilde\theta_n] = [\ombar_1:\dots:\ombar_n].$
Thus the functions $\gamma_{ij}$ define $\orbit$ on $\lag:$ we have
$\{I=\Ibar,\ \gamma_{ij} = 0\ \forall i,j\} = \orbit;$ indeed, the
vanishing of each $\gamma_{i,i+1}$ and of $\gamma_{n,1}$ suffices to define
$\orbit,$ and these $n$ functions may be taken as coordinates on $\lag$ in
a neighborhood of $\orbit.$ The central point of our argument will be that
the $\gamma_{ij}$ are ``propagating variables'' with derivatives along the
flow that, taken together, will suffice to give Lagrangian regularity.

Since the $\gamma_{kl}$ define $\orbit$ and are smooth on the annular
neighborhood $U$ where we have assumed regularity, there is a smooth cutoff
function $$\psi
:=\psi(\gamma_{12,},\gamma_{23},\dots,\gamma_{n-1,n},\gamma_{n,1})$$ with
$\psi=1$ on $\orbit$ and $\nabla \psi$ having its support on $\lag$
contained in $U.$ We may also arrange for $\psi$ to be the square of a
smooth function.  Let $\phi_\ep$ be a cutoff supported in $[-\ep, \ep],$
with smooth square root.

Let
\begin{equation}\label{adef}
a_{ij}(x) =   \psi \cdot \phi_\ep(\abs{I-\Ibar})\cdot \gamma_{ij}\cdot (\omega_i \ombar_j -\omega_j \ombar_i)
\end{equation}

We compute first that, where $\gamma_{ij} \in \CI,$
\begin{equation}\label{pgamma}
\{p,\gamma_{ij}\} = (\omega_i\ombar_j-\omega_j\ombar_i)
\end{equation}
(since $\gamma_{ij}$ is locally given by expressions of the form
$((\theta_i+2\pi k) \ombar_j-(\theta_j+2\pi l) \ombar_i)$ with $k,l$
\emph{fixed}) and hence that
$$
\{p,a_{ij}\}= \{p,\psi\} \phi_\ep (\abs{I-\Ibar}) \gamma_{ij} (\omega_i \ombar_j -\omega_j \ombar_i)+
\psi \phi_\ep (\abs{I-\Ibar}) (\omega_i\ombar_j-\omega_j\ombar_i)^2.
$$
We further note that as $\psi$ is a function of the $\gamma_{ij}$'s, by
\eqref{pgamma} the first term in this expression is a sum of
terms divisible by $$(\omega_{k_1} \ombar_{l_1} -\omega_{l_1} \ombar_{k_1})(\omega_{k_2} \ombar_{l_2} -\omega_{l_2} \ombar_{k_2})$$ for various
$k_i,l_i.$  Thus we may write
\begin{equation}\label{bracket}
\{p,a_{ij}\}= \sum e_k f_l+
\psi \phi_\ep (\abs{I-\Ibar}) (\omega_i\ombar_j-\omega_j\ombar_i)^2,
\end{equation}
where each $e_k$ and $f_l$ vanishes on $\lag$ and with support
intersecting $\lag$ only in $U.$

We will also employ a symbol that is \emph{invariant} under the flow:
for each $j=1,\dots, n,$ set
$$
w_j = \phi_\ep(\abs{I-\Ibar}) I_j.
$$

\section{Nondegeneracy}
Using a positive commutator argument, we will find that we can control
operators whose symbols are multiples of
$(\omega_i\ombar_j-\omega_j\ombar_i).$ These quantities vanish on $\lag,$
but our nondegeneracy hypothesis permits us to use them to control
Lagrangian regularity on $\lag.$  To see this, rewrite
$$
(\omega_i\ombar_j-\omega_j\ombar_i) = (\omega_i-\ombar_i)\ombar_j-(\omega_j
-\ombar_j)\ombar_i$$
and expand about $\lag$ in the $I$ variables, to rewrite this as
$$
\sum_k (\ombar_{ik}\ombar_j-\ombar_{jk}
\ombar_i) (I_k-\Ibar_k)+ O((I-\Ibar)^2)
$$

We now prove a key algebraic lemma:
\begin{lemma}\label{lemma:nondeg}
Let $v_1,\dots, v_n,$ and $v_{ij},$ $i,j=1,\dots,n$ be real numbers, with
$v_{ij}=v_{ji}.$ The functionals $\alpha_{ij}(x) = \sum_k(v_{ik}v_j-v_{jk}
v_i)x_k$ (for $i,j=1,\dots,n$) together with the covector $(v_1,\dots,v_n)$
span $(R^n)^*$ iff the matrix
\begin{equation}\label{thematrix}
\begin{pmatrix}
v_{11}& \dots & v_{1n} & v_1 \\
\vdots & \ddots & \vdots & \vdots\\
v_{n1} & \dots & v_{nn} & v_n\\
v_{1} & \dots & v_{n} & 0\\
\end{pmatrix}
\end{equation}
is nondegenerate.
\end{lemma}
\begin{proof}[Proof of Lemma]
We may assume that not all of the $v_i$'s are zero, as the result is
trivial in that case.

Let
$$
\vec\zeta_{ij} = \begin{pmatrix}
v_{1i}v_j-v_{1j}v_i\\
\vdots\\
v_{ni}v_j-v_{nj}v_i
\end{pmatrix}.
$$
Letting $A$ be the matrix with entries $v_{ij}$ and
$$
\vec u_{ij} = v_j e_i-v_i e_j
$$
where $e_i$ is the standard basis for $\RR^n,$
we have
$$
\vec\zeta_{ij} = A \vec u_{ij}.
$$
Let $U$ denote the span of the $\vec u_{ij}$'s.  Thus,
$$
U^\perp =  \bigcap_{i,j} \vec u_{ij}^\perp = \bigcap_{i,j} \RR \cdot \{\vec{w}
\in \RR^n| [w_i:w_j]=[v_i:v_j]\ \forall i,j\} = \RR \vec v
$$ where $\vec v = (v_1,\dots v_n)^t.$ Thus, $U = \vec v^\perp.$ Hence the
span of the $\vec \zeta_{ij}$ is of $A (\vec v^\perp).$ The assertion of
the lemma is then that $A(\vec v^\perp)$ and $\vec v$ are complementary iff
the matrix \eqref{thematrix} is nondegenerate.  This equivalence follows
from the observation that
$$
\begin{pmatrix}
A & \vec v \\
\vec v^t & 0 \\
\end{pmatrix} \cdot
\begin{pmatrix}
\vec w \\ z
\end{pmatrix}
= 
\begin{pmatrix}
A \vec w +z \vec v \\ \ang{\vec v,\vec w},
\end{pmatrix}
$$ hence \eqref{thematrix} has nontrivial nullspace iff there exists a
nonzero $\vec w \in \vec v^\perp$ with $A \vec w \in \RR \vec v.$
\end{proof}

\section{Proof of Theorem B}

We note, first of all, that in the special case when $n=2,$ a neighborhood
of any closed orbit $\orbit'\neq \orbit$ is itself an annular neighborhood
of $\orbit.$  Hence the special result for $n=2$ follows directly from the
general one.

We now prove Theorem B by induction on $k;$ we suppose it true for $k\leq
K-1$ (and note that for $k=0$ it is vacuous).

Let $A_{ij} \in \Psih(M)$ be self-adjoint, with symbol $a_{ij}$ constructed
above and vanishing subprincipal symbol.  Then we have by \eqref{bracket},
\begin{equation}\label{comm1} ih^{-3} [P, A_{ij}] = h^{-2} B_{ij}^2 + \sum_{k,l}
h^{-2}E_k F_l + R
\end{equation} with $B_{ij}$ self-adjoint with vanishing
subprincipal symbols, and 
\begin{equation}\label{bsymbol}
\sigma(B_{ij})=b_{ij} =({\psi}{\phi_\ep} (\abs{I-\Ibar}))^{1/2}\cdot
(\omega_i\ombar_j-\omega_j\ombar_i), 
\end{equation}
and with $E_k,$ $F_l$ characteristic on $\lag$ with the supports of
$\sigma(E_k),$ $\sigma(F_l)$ intersecting $\lag$ only in $U.$ ($R,$ $E_k,$
and $F_l$ of course depend on $i,j$ but we suppress these extra indices.)

Let $W_j$ have symbol $w_j$ constructed above, and be self-adjoint with
vanishing subprincipal symbol.  Then
$$
i h^{-3}[P, W_j] \in \Psi_h(M).
$$
For a multi-index $\alpha$ with $\abs{\alpha} = K-1,$
set
$$
Q_{ij}=A_{ij} W_1^{2\alpha_1}\dots W_n^{2\alpha_n},
$$
We will also need the operator denoted in multi-index notation
$$
W^\alpha = W_1^{\alpha_1}\dots W_n^{\alpha_{K-1}}
$$

Now we examine
\begin{equation}
i h^{-2K-1} \ang{ (P^* Q_{ij}-Q_{ij}P  u, u} = i
  h^{-2K-1} (\ang{Q_{ij}u,f}-\ang{f,Q_{ij}^* u}).
\end{equation}
For any $\delta>0,$ we may estimate the RHS by
$$
C_\delta\norm{h^{-K-1} f}^2+\delta(\norm{h^{-K} Q_{ij}u}^2+\norm{h^{-K} Q_{ij}^*u}^2).
$$ Note that both $Q_{ij}$ and $Q_{ij}^*$ are $(2K+1)$-fold products of
operators vanishing on $M,$ and that each contains the factors $A_{ij}$ and
$W^\alpha.$ By \eqref{adef} and \eqref{bsymbol}, $\sigma(A_{ij})$ is
divisible by $\sigma(B_{ij});$ thus, by elliptic regularity we may estimate
the RHS by
$$
C_\delta\norm{h^{-K-1} f}^2+C \delta\sum_{\alpha} \norm{h^{-K} B_{ij} W^\alpha u}^2+
\sum_{j=0}^{2K-2} h^{-j} \ang{D_j u,u}
$$ where $C$ is independent of $\delta,$ and each $D_j$ is a sum of
products of $j$ elements of $\md,$ all microsupported on $\supp A_{ij};$
these arise from commutator terms in which we have reordered products of
elements of $\md.$

Now we recall that $P^*-P \in h^2 \Psi_h(M),$ hence, by \eqref{comm1} we
may write
\begin{multline}\label{LHS}
\ang{ih^{-2K-1} (P^* Q_{ij}-Q_{ij}P)u,u} \\= \norm{h^{-K} B_{ij} W^\alpha
  u}^2  
+ \sum_{k,l}\ang{h^{-K} E_k W^\alpha u, h^{-K} F_l W^\alpha u}
+\sum_{j=0}^{2K-2} h^{-j} \ang{\tilde{D}_j u,u}
\end{multline}
with the $\tilde{D}_j$ sharing the properties of the $D_j$ above.

Putting together the information from our commutator, we now have, for all $\delta>0,$
\begin{multline}\label{ijestimate}
(1-C\delta)\norm{h^{-2K} B_{ij} W^\alpha u}^2 \\ \leq C_\delta\norm{h^{-K-1} f}^2+
\sum_{k,l}(\norm{h^{-K} E_k W^\alpha u}^2+\norm{h^{-K} F_l W^\alpha u}^2)
+\sum_{j=0}^{2K-2} h^{-j} \abs{\ang{\tilde{\tilde{D}}_j u,u}},
\end{multline}
with the $\tilde{\tilde{D}}_j$'s satisfying the same properties as $D_j$
above.  Each of the $E_k$ and $F_l$ terms is controlled by our hypothesis
on $\Slag^k u,$ while the $\tilde{\tilde{D}}_j$ terms are bounded by the
inductive assumption.

Now we use our nondegeneracy hypothesis as reflected in
Lemma~\ref{lemma:nondeg}.  Recall that $\lag \subset \Sigma,$ hence the
operator $P$ is characteristic on $\lag;$ moreover, we have $dp|_{\lag} =
\sum \ombar_k dI_k,$ hence Lemma~\ref{lemma:nondeg} tells us that $P$ and
$B_{ij},$ for $i,j=1,\dots,n,$ are a collection of operators fulfilling the
hypotheses of Proposition~\ref{prop:basis}.  Thus, adding together
equations \eqref{ijestimate} for all possible values of $i,j,$ and
multi-index $\alpha,$ together with terms involving $P$ rather than
$B_{ij}$ (which vanish up to commutators of $P$ with $W$'s), we obtain the
desired estimate, by Proposition~\ref{prop:basis}.\qed

\end{document}